\documentclass[preprint,10pt]{article}
\usepackage{cite}
\usepackage{url}
\usepackage{ifthen}
\usepackage{multicol}
\usepackage{amsmath,amssymb,amsthm}
\usepackage{latexsym,amssymb,amsmath}
\usepackage[utf8]{inputenc}
\usepackage{graphicx}
\usepackage{epstopdf}
\usepackage{subfig}
\usepackage{float}
\usepackage{color}
\newtheorem{theorem}{Theorem}
\newtheorem{corollary}{Corollary}

\newtheorem{definition}{Definition}

\newtheorem{remark}{Remark}

\numberwithin{equation}{section}
\newboolean{publ}

\newcommand{\ncm}{\newcommand}
\ncm{\maint}{\rlap{$\,\hspace{1.5790pt}=$}\int}
\ncm{\mcint}{\rlap{$\,\hspace{1.5790pt}-$}\int}

\title{Boundedness of normalization generalized differential operator of fractional formal }
\author{  Zainab E. Abdulnaby$^{1}$,  Rabha W. Ibrahim$^{2}$ and  Adem Kilicman$^{1,*}$   \\
          \small $^1$ Department of Mathematics, Universiti Putra Malaysia, \\
         \small  43400 UPM Serdang, Selangor, Malaysia\\
          \small $^2$ Faculty of Computer Science and Information Technology\\
          \small University Malaya, 50603, Malaysia
         \\
                       \small $^*$ corresponding 
}
\begin{document}
\maketitle

\begin{abstract} Many authors have considered and investigated generalized fractional differential operators.
 The main object of this present paper is to define a new generalized fractional differential operator $\mathfrak{T}^
 {\beta,\tau,\gamma},$ which generalized the Srivastava-Owa
 operators. Moreover, we investigate of the geometric properties such as univalency, starlikeness, convexity  for their
 normalization. Further, boundedness and compactness in some well known spaces, such as Bloch
  space for last mention operator also are considered. Our tool is
  based on the generalized hypergeometric function.
 \end{abstract}

\section{Introduction}
The study of fractional operators ( integral and differential) plays a vital and essential role in mathematical applied and
 mathematical analysis. To define a generalized  fractional differential operators and study their properties is one of
  important areas of current ongoing research in the geometric function theory and its concerning fields.
Many authors generalized fractional differential operators on previously known classes of analytic and univalent functions
to discover and produce new classes and to investigate various interesting properties of new classes, for example,
 (see \cite{z1,z18,z15,z24,z21}). In addition, several interesting applications of special function based definition
  for generalized fractional differential operators can be found in \cite{z10,z30,z31,z32}.

\medskip \noindent  let $ \mathcal{A}$ denote the class of functions  $ f(z)$  of the  form:
\begin{align}\label{de: ana2}
 f(z)= &  z+ \sum_{\kappa= 2}^{\infty} a_{\kappa} z^{\kappa} \end{align}
\noindent which are analytic functions in the open unit disk $ \mathbb{U} := \left\lbrace z \in \mathbb{C}: |z|<1 \right\rbrace $
and satisfy the condition $ f(0)=f^{\prime}(0)-1=0$. The subclass of $ \mathcal{A}$, consisting of all univalent(or ono-to-one)
functions $ f(z)$ in $ \mathbb{U}$ is denoted by $ \mathcal{S}$.
A function  $f(z) \in  \mathcal{A}$ is a  starlike functions  of order $ \lambda\,(0\leq \lambda<1)$, if and only if
 \begin{eqnarray}
 \Re \lbrace\frac{z f^{\prime}(z)}{f(z)}\rbrace  > \lambda,\quad   z \in \mathbb{U}
 \end{eqnarray}
   and denoted by
  $ \mathcal{S}_{\lambda}^{*}(m)$.
  Similarly, if $ f(z) \in \mathcal{A}$ and satisfies the following inequality
   \begin{eqnarray}
   \Re \lbrace \frac{z f^{\prime \prime}(z)}{ f^{\prime}(z)} +1\rbrace  > \lambda, \quad    z \in \mathbb{U}
   \end{eqnarray}
for some $ \lambda\,(0\leq \lambda<1)$, then $ f$ is said to be convex function of order $ \lambda$ and we denote this class
 by $ \mathcal{K}_{\lambda}(m)$.
\begin{theorem} \label{remark}\rm{(Bieberbach's Conjecture \cite{z2}) If the function $f (z)$ defined by \eqref{de: ana2}
 is in the class $\mathcal{S}^{*}$ then $ | a_{\kappa}| \leq \kappa $  for all  $ \kappa \geq 2 $  and if it is  in the
  class $ \mathcal{K}$  then $| a_{\kappa}|\leq 1 $ for all $ \kappa \geq 2$.}\
\end{theorem}
\noindent For the function  $f(z)$ defined by \eqref{de: ana2} and  $ h(z)= z + \sum_{\kappa= 2}^{\infty} b_{\kappa} z^{\kappa},$
 the convolution (or Hadamard product) $  f* h $ is given by
\begin{align}\label{de: Conv}
 f * h(z) =  z + \sum_{\kappa= 2}^{\infty} a_{\kappa} b_{\kappa} z^{\kappa}
 \end{align}
\\
\medskip \noindent The  operator  $\mathcal{O}_{z}^{\beta,\tau}$ is defined in terms of Riemann-Liouville fractional differential
 operator $ \mathcal{D}_{z}^{\beta-\tau}$ as
\begin{align}
\mathcal{O}_{z}^{\beta,\tau}f(z)=\frac{\Gamma(\tau)}{\Gamma(\beta)}z^{1-\tau} \mathcal{D}_{z}^{\beta-\tau} z^{\beta-1} f(z)
\quad ( z \in \mathbb{U}),
\end{align}
This operator given by Tremblay \cite{z3}. Recently, Ibrahim \cite{z4} extended Tremblay's operator in terms of Srivastava-Owa
fractional differential of $ f(z)$ of order $(\beta-\tau)$ and defined as follows
\begin{align}\label{Tre}
\mathfrak{T}_{z}^{\beta,\tau}f(z)=\frac{\Gamma(\tau)}{\Gamma(\beta)}z^{1-\tau} D_{z}^{\beta-\tau} z^{\beta-1} f(z),
\end{align}
or, equivalent
\begin{align}
\mathfrak{T}_{z}^{\beta,\tau}f(z)=\frac{\Gamma(\tau)}{\Gamma(\beta)}z^{1-\tau} \frac{\rm{d}}{\rm{d}z} \int_{0}^{z}(z-\zeta)
^{\tau-\beta} \zeta^{\tau-1} d\zeta  \quad ( z \in \mathbb{U}).
\end{align}
Often, the generalized fractional differential operators and their applications associated with special functions
 ( see \cite{z22}). So, the Fox-Wright $ {}_{p}\Psi_{q}$ generalization of the hypergeometric $ {}_{p}F_{q}$  function is considered one of the important special function in geometric function theory defined by
  \begin{eqnarray}\label{Fox}
{}_{p}\Psi_{q}[z]= {}_{p}\Psi_{q}
 \begin{bmatrix}
  (a_{j},A_{j})_{1,p}; \\[0.1em]                                 \qquad \qquad z \\[0.1em]
 (b_{j},B_{j})_{1,q};
           \end{bmatrix}
     &=& \sum_{\kappa=0}^{\infty} \frac{\Pi_{j=1}^{p} \Gamma(a_{j} +\kappa A_{j})}{\Pi_{j=1}^{q} \Gamma(b_{j}+\kappa B_{j})(1)_{\kappa}}z^{\kappa}.
 \end{eqnarray}
 where $ a_{j}, b_{j} $ are parameters in complex plane $\mathbb{C}$. $ A_{j}>0 $, $ B_{j}>0 $ for all $ j=1, \ldots, q $ and  $ j=1,\ldots,p $, such that   $ 0\leq 1 + \sum_{j=1}^{q} B_{j} - \sum_{i=1}^{p} A_{i} $  for fitting values $|z| <1 $
and it is well know that
\[ {}_{p}\Psi_{q} \begin{bmatrix}
(a_{j},1)_{1,p};\\[0.1em]
\qquad \quad z\\[0.1em]
(b_{j},1)_{1,q}
\end{bmatrix} = \Delta^{-1} {}_{p}F_{q}(a_{j},b_{j};z)
\]
 where
 \[ \Delta= \frac{\Pi_{j=1}^{q} \Gamma(b_{j})}{\Pi_{j=1}^{p} \Gamma(a_{j})}\]
 and $ {}_{p}F_{q}$ is the generalized hypergeometric function. The \textit{Pochhammer} symbol denoted by  $(\rho)_{\kappa}$ and defined as follows
\begin{eqnarray}\label{q:p}
(\rho)_{\kappa} :=\dfrac{\Gamma(\rho + \kappa)}{\Gamma(\rho)} =
\begin{cases}
\rho (\rho +1)...(\rho + \kappa -1)& (\kappa \in \mathbb{N}; \rho  \in \mathbb{C}),\\
1 &(\kappa = 0 ; \rho \in \mathbb{C}\setminus \left\lbrace {0}\right \rbrace ),
\end{cases}
\end{eqnarray}
where $ \Gamma$ is the well know Gamma function.\\
  \noindent
   In the present paper, the generalized Tremblay operator $ \mathfrak{T}_{z}^{\beta ,\tau, \gamma}$ of analytic function is defined. Also, the univalence properties of normalization generalized operator  are investigated and proved. Further, the boundedness and compactness of the last operator are studied.
 \section{Background and Results}

\noindent  In this section, we consider the generalized Tremblay type fractional differential operator and then we determine the generalized fractional differential of some special functions. For this main purpose, we begin by recalling the  Srivstava-Owa fractional differential operators of $ f(z)$ of order $\beta$ defined by
 \begin{eqnarray}
D^{\beta}_{z} f(z) := \frac{1}{\Gamma(1-\beta)}\frac{{\rm{d}}}{{\rm{d}}z}\int_{0}^{z} f(\zeta) (z- \zeta)^{-\beta} d\zeta,
\end{eqnarray}
 where $ 0 \leq \beta < 1 $, and the function $ f(z)$ is analytic in simply-connected region of the complex $z$-plane containing the origin and the multiplicity of $ (z-\zeta)^{-\beta}$ is removed by requiring $ \log(z-\zeta)$ to be real when $ (z-\zeta)>0$ (see \cite{z5,z6}). Then under the conditions of the above definition the Srivstava-Owa fractional differential of $ f(z)= z^{\kappa}$ is defined by
  \begin{align}\label{1}
  D_{z}^{\beta}\lbrace z^{\kappa} \rbrace= \frac{\Gamma(\kappa+1)}{\Gamma(\kappa-\rho+1)}z^{\kappa-\beta}.
  \end{align}
The theory of fractional integral and differential operators has found significant importance applications in various areas, for example \cite{z9}. Recently, many mathematicians have developed various generalized fractional differential of Srivastava-Owa type, for example, see \cite{z7,z8}. Further, we consider a new generalized Srivastava-Owa type fractional differential formulas which is recently appeared.
\begin{definition}\rm{\cite{z1}
The generalized Srivastava-Owa fractional differential  of $ f(z)$ of order $ \beta$ is defined by
\begin{align}\label{3}
\mathcal{D}_{z}^{\beta,\gamma} f(z):= \frac{(\gamma +1)^{\beta}}{\Gamma(1- \beta)}\frac{{\rm{d}}}{{\rm{d}}z}\int_{0}^{z} (z^{\gamma+1} - \zeta^{\gamma+1})^{-\beta} \zeta^{\gamma} f(\zeta) d\zeta,
\end{align}
 where $0 \leq \beta <1$, $\gamma > 0$ and $ f(z)$ is analytic in simply-connected region of the complex z-plane $ C $ containing the origin, and the multiplicity of $(z^{\gamma+1} - \zeta^{\gamma+1})^{-\beta}$ is removed by requiring $ \log(z-\zeta)$ to be real when $(z-\zeta)>0$. In particular, the generalized Srivastava-Owa fractional differential of function $ f(z)= z^{\kappa}$ is defined by
 \begin{align*} \label{4}
\mathcal{ D}^{\beta,\gamma}_{z} \left\lbrace z^{\kappa}\right\rbrace  = \frac{(\gamma+1)^{\beta-1}\Gamma \left(\frac{\kappa}{\gamma+1}+1\right)}{\Gamma\left(\frac{\kappa}{\gamma+1}+1-\beta\right)} z^{(1-\beta)(\gamma+1)+\kappa-1}.
 \end{align*}
    }\end{definition}
\noindent Now, we present the generalized operator of type fractional differential as follows:
  \begin{definition}\label{GT4}\rm{
The generalized fractional differential of $f(z)$ of two parameters $ \beta$ and $\tau$ is defined by
 \begin{equation}\label{AT}
  \mathfrak{T}_{z}^{\beta ,\tau, \gamma}f(z) :=\frac{(\gamma +1)^{\beta-\tau} \Gamma(\tau)}{\Gamma(\beta)\Gamma(1-\beta-\tau)} \bigg( z^{1-\tau} \frac{{\rm{d}}}{{\rm{d}}z} \bigg)  \int_{0}^{z} \frac{\zeta^{\gamma+\beta-1}f(\zeta) }{(z^{\gamma +1} - \zeta^{\gamma +1})^{\beta -\tau}}\,d\zeta,
 \end{equation}
 $$(\gamma\geq 0;\,0< \beta \leq 1;\, 0< \tau \leq 1;\, 0\leq \beta-\tau <1),$$
where the function $ f(z)$ is analytic in simple-connected region of the complex z-plane $ \mathbb{C}$ containing the origin, and the multiplicity of $(z^{\gamma+1}-\zeta^{\gamma+1})^{-\beta+\tau} $ is removed by requiring $\log (z^{\gamma+1} - \zeta^{\gamma+1})$ to be non-negative when $( z^{\gamma+1}-\zeta^{\gamma+1}) > 0 $. } \end{definition}
\begin{remark}\rm{ For $ f(z) \in \mathcal{A}$, we have
\begin{itemize}
\item[i-] When $ \gamma=0$ in \eqref{AT}, is reduce to the classical known one \eqref{Tre}.
\item[ii-] When $ \tau =\beta$, then
 \begin{align}
\mathfrak{T}_{z}^{\beta,\beta,\gamma} f(z)=  f(z).
\end{align}
\end{itemize}
}\end{remark}

\noindent We investigate the generalized fractional differential of the function $f(z)= z^{\nu};\,\nu \geq 0$.
    \begin{theorem}\label{e}\rm{
 Let $ 0 \leq \beta-\tau<1 $ for some $ 0< \beta \leq 1$; $ 0 < \tau\leq 1$ and  $\upsilon \in N $, then we have
 \begin{eqnarray}\label{zp}
 \mathfrak{T}_{z}^{\beta,\tau,\gamma} \lbrace z^{\upsilon}\rbrace = \frac{(\gamma+1)^{\beta-\tau} \Gamma(\frac{\upsilon+\beta-1}{\gamma+1}+1)\Gamma(\tau)}{\Gamma(\frac{\upsilon+\beta-1}{\gamma +1}+1 -\beta+\tau)\Gamma(\beta)} z^{(1-\beta+ \tau)\gamma+\upsilon}.
\end{eqnarray}
   \begin{proof}
   Applying \eqref{AT} in Definition \ref{GT4} to the function $ z^{\upsilon}$, we obtain
\begin{small}   \begin{align}
   \mathfrak{T}_{z}^{\beta,\tau,\gamma} \lbrace z^{\upsilon}\rbrace &= \frac{(\gamma +1)^{\beta-\tau} \Gamma(\tau)}{\Gamma(\beta)\Gamma(1-\beta+\tau)} \bigg( z^{1-\tau} \frac{{\rm{d}}}{{\rm{d}}z} \bigg)  \int_{0}^{z} {\zeta^{\gamma+\beta+ \upsilon-1} }{(z^{\gamma +1} - \zeta^{\gamma +1})^{\tau -\beta}}  d\zeta ,\nonumber
   \\
     & = \frac{(\gamma +1)^{\beta-\tau} \Gamma(\tau)}{\Gamma(\beta)\Gamma(1-\beta+\tau)} \bigg( z^{1-\tau} \frac{{\rm{d}}}{{\rm{d}}z} \bigg)  \int_{0}^{z}z^{(\gamma +1)(\tau-\beta)} \zeta^{\gamma+\beta+ \upsilon-1} \left(1 - \frac{\zeta^{\gamma +1}}{z^{\gamma +1}}\right)^{\tau -\beta}  d\zeta ,\nonumber
   \end{align}\end{small}
   Let use the substitution $ w:=(\frac{\zeta}{z})^{\gamma+1}$ in this expression, then have
   \begin{small} \begin{align}
    \mathfrak{T}_{z}^{\beta,\tau,\gamma}\lbrace z^{\upsilon}\rbrace
     & = \frac{(\gamma +1)^{\beta-\tau-1} \Gamma(\tau)}{\Gamma(\beta)\Gamma(1-\beta+\tau)} \bigg(z^{1-\tau} \frac{{\rm{d}} }{{\rm{d}}z} z^{(\gamma+1)(\tau-\beta)+\gamma+\upsilon+\beta}  \bigg) \int_{0}^{1} w^{\frac{\beta+ \upsilon-1}{\gamma+1}} (1 - w)^{\tau -\beta} dw\nonumber
     \\
     & = \frac{[(\gamma+1)(\tau-\beta)+\gamma+\upsilon+\beta](\gamma +1)^{\beta-\tau-1} \Gamma(\tau)}{\Gamma(\beta)\Gamma(1-\beta+\tau)}\nonumber\\
     &\quad \qquad \qquad \qquad \qquad \times \bigg( z^{(1-\beta+ \tau)\gamma+\upsilon}  \bigg) B(\frac{ \upsilon+\beta-1}{\gamma+1}+1,1-\beta +\tau),\label{B}
     \\
     &= \frac{(\gamma+1)^{\beta-\tau} \Gamma(\frac{\upsilon+\beta-1}{\gamma+1}+1)\Gamma(\tau)}{\Gamma(\frac{\upsilon+\beta-1}{\gamma +1}+1-\beta+\tau)\Gamma(\beta)} z^{(1-\beta+ \tau)\gamma+\upsilon}.\nonumber
      \end{align}\end{small}
    where $ B(.,.)$ in \eqref{B} is the  Beta function defined by [ ]
    \[
    B(u,\upsilon)= \begin{cases}  \int_{0}^{1} \eta^{u}(1-\eta)^{\upsilon-1} d\eta \qquad ( u >0;\, \upsilon>0), \\[0.6em]
\frac{\Gamma(u)\Gamma(\upsilon)}{\Gamma(u +\upsilon)}\qquad \qquad\,\,\, \qquad (u,\,\upsilon \in \mathbb{C} \setminus \mathbb{Z}_{0}^{-}).
 \end{cases}
    \] Thus, the proof Theorem \ref{e} is completed.
              \end{proof}
 }\ \end{theorem}
 \subsection{Applications by using special functions}
 In this subsection we replace the normalized of function $ f(z)$ by some special functions in Theorem \ref{e}. Begin let recall  generalized of \textit{Koebe} function $ f(z)=z /(z-1)^{-\alpha}$ as follows:
  \begin{theorem}\label{The Koby}\rm{
 Let $ f(z)= z(1-z)^{-\alpha}$, $ \alpha \geq 1$ and $ z \in \mathbb{U}$, then
 \begin{eqnarray}
  \mathfrak{T}_{z}^{\beta,\tau,\gamma} f(z)= \frac{(\gamma+1)^{\beta-\tau}\Gamma{(\tau)}}{\Gamma(\beta)}z^{(1-\beta+ \tau)\gamma+1}
 {}_{3}\Psi_{2} \begin{scriptsize}
 \begin{bmatrix}
  (\alpha+1,1),(1,1), (1+ \frac{\beta}{\gamma+1 },\frac{1}{\gamma +1});\\
\qquad   \qquad  \qquad \qquad \qquad \qquad \qquad z \\
     (2,1), (1-\beta +\tau + \frac{\beta}{\gamma +1} ,\frac{1}{\gamma +1});
           \end{bmatrix}
\end{scriptsize}.\nonumber
\end{eqnarray}
\begin{proof} Since
\[ \mathfrak{T}_{z}^{\beta,\tau,\gamma} \lbrace z(1-z)^{-\alpha}\rbrace =
   \mathfrak{T}_{z}^{\beta,\tau,\gamma} \lbrace\sum_{\kappa=1}^{\infty} \frac{(\alpha)_{\kappa}}{(1)_{\kappa}}z^{\kappa}\rbrace\] then, by using Theorem \ref{e} we obtain
\begin{eqnarray}
&=&\frac{(\gamma +1)^{\beta-\tau} \Gamma(\tau)}{\Gamma(\beta)\Gamma(1-\beta+\tau)} \bigg( z^{1-\tau} \frac{{\rm{d}}}{{\rm{d}}z} \bigg)  \int_{0}^{z} {\zeta^{\gamma+\beta-1} }{(z^{\gamma +1} - \zeta^{\gamma +1})^{\tau -\beta}} \sum_{\kappa=1}^{\infty} \frac{(\alpha)_{\kappa}}{(1)_{\kappa}}\zeta^{\kappa} d\zeta ,\nonumber\\[0.6em]
 &=& \frac{(\gamma+1)^{\beta-\tau}\Gamma{(\tau)}}{\Gamma(\beta)}z^{(1-\beta+ \tau)\gamma}\sum_{\kappa=1}^{\infty}  \frac{\Gamma(\kappa+\alpha)\Gamma(\frac{\kappa}{\gamma+1}+ \frac{\beta}{\gamma+1}-\frac{1}{\gamma+1}+1)}{\Gamma(\alpha)\Gamma(\frac{\kappa}{\gamma +1}+ \frac{\beta}{\gamma+1}-\frac{1}{\gamma+1}+1-\beta+\tau)(1)_{\kappa}}  z^{\kappa}  \nonumber\\[0.6em]
 &=& \frac{\alpha(\gamma+1)^{\beta-\tau}\Gamma{(\tau)}}{\Gamma(\beta)}z^{(1-\beta+ \tau)\gamma+1}\sum_{\kappa=0}^{\infty}  \frac{\Gamma(\kappa+\alpha+1)\Gamma(\kappa+1)\Gamma(\frac{\kappa}{\gamma+1}+ \frac{\beta}{\gamma+1}+1)}{\Gamma(\alpha+1)\Gamma(\frac{\kappa}{\gamma +1}+ \frac{\beta}{\gamma+1}+1-\beta+\tau)\Gamma(\kappa+2)} \frac{ z^{\kappa} }{(1)_{\kappa}} \nonumber\\[0.6em]
  &=& \frac{\alpha(\gamma+1)^{\beta-\tau}\Gamma{(\tau)}}{\Gamma(\beta)}z^{(1-\beta+ \tau)\gamma+1}
 {}_{3}\Psi_{2} \begin{scriptsize}
 \begin{bmatrix}
  (\alpha+1,1),(1,1), (1+ \frac{\beta}{\gamma+1 },\frac{1}{\gamma +1});\\
\qquad   \qquad  \qquad \qquad \qquad \qquad \qquad z \\
     (2,1), (1-\beta +\tau + \frac{\beta}{\gamma +1} ,\frac{1}{\gamma +1});
           \end{bmatrix}
\end{scriptsize}.\nonumber
 \end{eqnarray}
\end{proof}
      } \end{theorem}
 Setting $ \alpha =2$ in Theorem \ref{The Koby}
 , we obtain
  \begin{corollary}\rm{
 Let $ f(z) = z(1-z)^{-2}$, then
 \begin{eqnarray}
  \mathfrak{T}_{z}^{\beta,\tau,\gamma} f(z)= \frac{2(\gamma+1)^{\beta-\tau}\Gamma{(\tau)}}{\Gamma(\beta)}z^{(1-\beta+ \tau)\gamma+1}
 {}_{3}\Psi_{2} \begin{scriptsize}
 \begin{bmatrix}
  (3,1),(1,1), (1+ \frac{\beta}{\gamma+1 },\frac{1}{\gamma +1});\\
\qquad   \qquad  \qquad \qquad \qquad \qquad \qquad z \\
     (2,1), (1-\beta +\tau + \frac{\beta}{\gamma +1} ,\frac{1}{\gamma +1});
           \end{bmatrix}
\end{scriptsize}.\nonumber
\end{eqnarray}
 } \end{corollary}
 Setting $ \alpha=1$, in Theorem \ref{The Koby}, we obtain
  \begin{corollary}\rm{
 Let $ f(z)= z(1-z)^{-1}$ then we have
   \begin{eqnarray}
  \mathfrak{T}_{z}^{\beta,\tau,\gamma} f(z)= \frac{(\gamma+1)^{\beta-\tau}\Gamma{(\tau)}}{\Gamma(\beta)}z^{(1-\beta+ \tau)\gamma+1}
 {}_{3}\Psi_{2} \begin{scriptsize}
 \begin{bmatrix}
  (2,1),(1,1), (1+ \frac{\beta}{\gamma+1 },\frac{1}{\gamma +1});\\
\qquad   \qquad  \qquad \qquad \qquad \qquad \qquad z \\
     (2,1), (1-\beta +\tau + \frac{\beta}{\gamma +1} ,\frac{1}{\gamma +1});
           \end{bmatrix}
\end{scriptsize}.\nonumber
\end{eqnarray}
 } \end{corollary}
 \bigskip \noindent The following Theorems is seen to immediately follow from Theorems \ref{e} and \ref{The Koby}.

  \begin{theorem}  \rm{ Let $ 0 < \beta\leq 1$ and $ 0 < \tau\leq 1$ such that  $0 \leq \beta-\tau < 1$. Then we have
 \begin{eqnarray*}
 \mathfrak{T}_{z}^{\beta,\tau,\gamma}\lbrace z\,e^{z}\rbrace &=& \tfrac{(\gamma+1)^{\beta-\tau} \Gamma(\tau)}{\Gamma(\beta)\Gamma(1-\beta+\tau)} \bigg( z^{1-\tau} \frac{{\rm{d}}}{{\rm{d}}z} \bigg)  \int_{0}^{z} \zeta^{\gamma+\beta-1} (z^{\gamma+1} + \zeta^{\gamma+1})^{\tau-\beta} \sum_{\kappa=1}^{\infty}\frac{ (z\zeta)^{\kappa}}{(1)_{\kappa}} d\zeta,
 \\[0.6em]
 &=&  \frac{(\gamma+1)^{\beta-\tau} \Gamma(\tau)}{\Gamma(\beta)}   z^{(1-\beta+\tau)\gamma}  \sum_{\kappa=1}^{\infty} \frac{\Gamma(\frac{\kappa+\beta-1}{\gamma+1}+1)}{\Gamma(\frac{\kappa+\beta-1}{\gamma+1}+1-\beta+\tau)(1)_{\kappa}}z^{\kappa}.
 \\[0.6em]
  &=&  \frac{(\gamma+1)^{\beta-\tau} \Gamma(\tau)}{\Gamma(\beta)}   z^{(1-\beta+\tau)\gamma+1}  \sum_{\kappa=0}^{\infty} \frac{\Gamma(\kappa+1)\Gamma(\frac{\kappa+\beta}{\gamma+1}+1)}{\Gamma(\frac{\kappa+\beta}{\gamma+1}+1-\beta+\tau)\Gamma(\kappa+2)}\frac{z^{\kappa}}{(1)_{\kappa}}.
 \\[0.6em]
 &=&  \frac{(\gamma+1)^{\beta-\tau} \Gamma(\tau)}{\Gamma(\beta)}   z^{(1-\beta+\tau)\gamma+1}
 {}_{2}\Psi_{2} \begin{scriptsize}
 \begin{bmatrix}
 (1,1), (1+ \frac{\beta}{\gamma+1 }, \frac{1}{\gamma +1}); \\
\qquad   \qquad  \qquad \qquad \qquad  \qquad z \\
   (2,1), (1-\beta +\tau + \frac{\beta}{\gamma +1} , \frac{1}{\gamma +1});            \end{bmatrix}
\end{scriptsize}.
\end{eqnarray*}
 }\end{theorem}
\bigskip
 \begin{theorem}\rm{
 Let $ 0 < \beta \leq  1$ and $0 < \tau \leq  1$ such that $ 0 \leq \beta-\tau < 1$. Then, we obtain
 \begin{eqnarray*}
 \mathfrak{T}_{z}^{\beta,\tau,\gamma}\lbrace z {}_{1}F_{1}(\alpha,\lambda;z)\rbrace = \frac{(\gamma+1)^{\beta-\tau} \Gamma(\tau)}{\Gamma(\beta)\Gamma(1-\beta+\tau)} z^{(1-\beta+\tau)\gamma} \sum_{\kappa=1}^{\infty}\tfrac{(\alpha)_{\kappa}}{(\lambda)_{\kappa}(1)_{\kappa}} B( \tfrac{\kappa+\beta-1}{\gamma+1}+1, 1-\beta+\tau) z^{\kappa}
\nonumber \\[0.6em]
=  \frac{\alpha (\gamma+1)^{\beta-\tau}\Gamma(\tau) }{\lambda\Gamma(1+\tau)} z^{(1-\beta+\tau)\gamma+1} \sum_{\kappa=0}^{\infty}\tfrac{(1)_{\kappa}(\alpha+1)_{\kappa}}{(\lambda+1)_{\kappa}(2)_{\kappa}} \frac{ B( \frac{\kappa+\beta}{\gamma+1}+1, 1-\beta+\tau)}{B(\beta,1-\beta +\tau)} \frac{z^{\kappa}}{(1)_{\kappa}},
 \end{eqnarray*}
 where ${}_{1}F_{1}(\alpha,\lambda;z)$ is  the confluent hypergeometric \textit{Kummer} function see \cite{z10} defined by 
 \begin{eqnarray*}
  {}_{1}F_{1}(\alpha,\lambda;z):&=& \frac{\Gamma(\lambda)}{\Gamma(\alpha)\Gamma(\lambda-\alpha)}\int_{0}^{z} t^{\alpha-1} (1-t)^{\lambda-\alpha-1} e^{zt}dt,\\[0.6em]
  &=&\sum_{\kappa=0}^{\infty}\frac{(\alpha)_{\kappa}}{(\lambda)_{\kappa}(1)_{\kappa}} z^{\kappa}.
\end{eqnarray*}

and $(\alpha)_{\kappa},\, (\lambda)_{\kappa}$ are given by \eqref{q:p}.
 }\end{theorem}
\begin{theorem}\rm{
 Let $ 0 < \beta \leq  1$ and $0 < \tau \leq  1$ such that $ 0 \leq \beta-\tau < 1$. Then, we obtain
 \begin{eqnarray*}
\mathfrak{T}_{z}^{\beta,\tau,\gamma}\lbrace z \Omega_{\alpha,\lambda,\rho}(z,s,r)\rbrace 
&=&\tfrac{(\gamma+1)^{\beta-\tau} \Gamma(\tau)}{\Gamma(\beta)\Gamma(1-\beta+\tau)} z^{(1-\beta+\tau)\gamma } \sum_{\kappa=1}^{\infty}\tfrac{(\alpha)_{\kappa}(\lambda)_{\kappa}}{(\rho)_{\kappa}(1)_{\kappa}} B( \tfrac{\kappa+\beta-1}{\gamma+1}+1, 1-\beta+\tau) \frac{z^{\kappa}}{(\kappa+a)^{s}}
\\[0.6em]
&=&\tfrac{\alpha \lambda(\gamma+1)^{\beta-\tau} \Gamma(\tau)}{\rho\Gamma(\beta)\Gamma(1-\beta+\tau)} z^{(1-\beta+\tau)\gamma+1 } \sum_{\kappa=0}^{\infty}\tfrac{(\alpha+1)_{\kappa}(\lambda+1)_{\kappa}}{(\rho+1)_{\kappa}\Gamma(\kappa+2)} B( \tfrac{\kappa+\beta}{\gamma+1}+1, 1-\beta+\tau) \tfrac{z^{\kappa}}{(\kappa + 1 +a)^{s}}
 \end{eqnarray*}

where extended general Hurwitz-Lerch Zeta function was introduced in \cite{z78} by
\begin{eqnarray*}
 \Omega_{\alpha,\lambda,\rho}(z,s,a) :&=& \sum_{\kappa=0}^{\infty}\frac{(\alpha)_{\kappa}(\lambda)_{\kappa}}{(\rho)_{\kappa}(1)_{\kappa}} \frac{z^{\kappa}}{(\kappa+a)^{s}}
\end{eqnarray*}
where $ \rho,\, a \in \mathbb{Z}\setminus\lbrace0,-1,-2,\cdots\rbrace,\, s \in \mathbb{C},\, \Re(s)>0$ when $|z|<1$,  for more details about this function see \cite{z77}.

 }\end{theorem}
 \section{Normalization generalized operator}
 In this section we normalize the generalized operator $\mathfrak{T}_{z}^{\beta,\tau,\gamma}$ of type fractional differential of analytic univalent functions  in $ \mathbb{U}$ and defined in two terms as follows:

\bigskip
\noindent  Let the following conditions to be realized:
 \begin{align}\label{cond}
  0\leq  \beta-\tau <1, \quad  \gamma \geq 0,\end{align}
\noindent we defined  the operator $\Theta^{\beta,\tau,\gamma}f(z): \mathcal{A}\rightarrow\mathcal{A}$ by \begin{eqnarray}\label{extend}
\Theta^{\beta,\tau,\gamma} f(z)&=& \frac{z^{(\beta-\tau-1)\gamma}\Gamma(\frac{\beta}{\gamma+1}+1-\beta+\tau)\Gamma(\beta)}{(\gamma+1)^{\beta-\tau}  \Gamma(\frac{\beta}{\gamma+1}+1)\Gamma(\tau)} \left\lbrace \mathfrak{T}_{z}^{\beta,\tau,\gamma} f(z) \right\rbrace \nonumber
 \\[0.6em]
&=&
   z+ \sum_{\kappa=2}^{\infty} \Phi_{\beta,\tau,\gamma}({\kappa}) \, a_{\kappa} z^{\kappa}
\end{eqnarray}
      where
   \begin{eqnarray}\label{Psi}
   \Phi_{\beta,\tau,\gamma}({\kappa}) := \frac{\Gamma(\frac{\beta}{\gamma+1}+1-\beta+\tau)\Gamma(\frac{\kappa+\beta-1}{\gamma+1}+1)}{ \Gamma(\frac{\beta}{\gamma+1}+1)\Gamma(\frac{\kappa+\beta-1} {\gamma+1}+1-\beta+\tau)}.
   \end{eqnarray}
 Not that \[ \Theta^{\beta,\tau,\gamma} f(0)=0\]
Next we employ the well known method of convolution product two functions of analytic univalent functions for $ \Theta^{\beta,\tau\gamma} $ and  define as follows:
 \begin{eqnarray}\label{q: G}
\Theta^{\beta,\tau,\gamma} f(z)
&=& z +
 \sum_{\kappa=2}^{\infty}\frac{\Gamma(\frac{\beta}{\gamma+1}+1-\beta+\tau)\Gamma(\frac{\kappa+\beta-1}{\gamma+1}+1)}{ \Gamma(\frac{\beta}{\gamma+1}+1)\Gamma(\frac{\kappa+\beta-1} {\gamma+1}+1-\beta+\tau)} a_{\kappa} z^{\kappa} \nonumber\\[0.6em]
&=& z \sum_{\kappa=0}^{\infty}\frac{\Gamma(\kappa+1)
 \Gamma(\frac{1+\beta}{\gamma+1}+1-\beta+\tau)\Gamma(\frac{\kappa+\beta}{\gamma+1}+1)}{ \Gamma(\frac{1+\beta}{\gamma+1}+1)\Gamma(\frac{\kappa+\beta} {\gamma+1}+1-\beta+\tau)(1)_{\kappa}(1)_{\kappa}} a_{\kappa}z^{\kappa}
 \nonumber
  \\[0.6em]
 &=& \frac{\Gamma(\frac{1+\beta}{\gamma+1}+1-\beta+\tau)}{\Gamma(\frac{1+\beta}{\gamma+1}+1)} z\, {}_{2}\Psi_{1}
 \begin{scriptsize}
 \begin{bmatrix}
  (1,1), (1+\frac{\beta}{\gamma+1}, \frac{1}{\gamma +1});
  \\
 \qquad \qquad \qquad \qquad \qquad z
 \\(1-\beta +\tau+\frac{\beta}{\gamma+1},\frac{1}{\gamma +1});       \end{bmatrix} * f(z).
\end{scriptsize}
 \end{eqnarray}
Therefore, \begin{eqnarray}\label{normalization}
\Theta^{\beta,\tau,\gamma} f(z):
= \frac{\Gamma(\frac{1+\beta}{\gamma+1}+1-\beta+\tau)}{\Gamma(\frac{1+\beta}{\gamma+1}+1)} z\,{}_{2}\Psi_{1} [z]* f(z)
\end{eqnarray}
  Now we investigate to study the univalence properties for operator  $\Theta^{\beta,\tau,\gamma}$
   \begin{theorem}\label{Theo 1}\rm{
  Let $ f\in \mathcal{S}$. If the following conditions satisfied
   \begin{itemize}
  \item[(i)] For $ 0 \leq \beta <0 $ , $ \tau >0 $ and $ \geq \beta-\tau <1$.
  \item[(ii)] $  \rho_{i} > 0,  i= 1, \ldots, p $ and $ \lambda_{j} > 0,  j= 1, \ldots, q $; $ p\leq q + 1$,
  \end{itemize}
   \noindent then the operator $\Theta ^{\beta,\tau,\gamma} f(z) \in \mathcal{S}$  in open unite disk $ \mathbb{U} $.
        \begin{eqnarray*}
{}_{2}\Psi_{1} \begin{scriptsize}
\begin{bmatrix}
  (3,1), (1+\frac{\beta}{\gamma+1 }+\frac{1}{\gamma+1} ,\frac{1}{\gamma +1});   \\
\qquad \qquad \qquad \qquad \qquad  \qquad \quad 1 \\
       (1-\beta +\tau+\frac{\beta}{\gamma +1} +\frac{1}{\gamma+1},\frac{1}{\gamma +1});  \\
           \end{bmatrix}\end{scriptsize} +
 {}_{2}\Psi_{1}
\begin{scriptsize} \begin{bmatrix}
  (2,1), (1+\frac{\beta}{\gamma+1 } ,\frac{1}{\gamma +1});            \\
 \qquad \qquad \qquad \qquad \qquad  1\\
       (1-\beta +\tau+\frac{\beta}{\gamma +1} ,\frac{1}{\beta +1});
           \end{bmatrix} \end{scriptsize}  < 2  \bigg(\tfrac{ \Gamma(\frac{\beta}{\gamma+1}+1)}{\Gamma(\frac{\beta}{\gamma+1}+1-\beta+\tau)} \bigg).                \end{eqnarray*}

\begin{proof}
Suppose the function $ f \in \mathcal{S}$ and let \[
 \Theta^{\beta,\tau,\gamma} f(z) =  z +       \sum_{\kappa=2}^{\infty} w_{\kappa}
  z^{\kappa}
\]
be defined by equality \eqref{extend}, where
 where
\[
  w_{\kappa}:=\Phi_{\beta,\tau,\gamma}(\kappa)\, a_{\kappa}
\]
   and the function $\Phi_{\beta,\tau,\gamma}(\kappa)$ is defined  by \eqref{Psi}. To prove that the operator $\Theta^{\beta,\tau,\gamma}$ preserves the class $ \mathcal{S}$ we require the following sufficient condition

\[\ell_{1};= \sum_{\kappa=2}^{\infty} \kappa\, |w_{\kappa}| = \sum_{\kappa=2}^{\infty}\, \kappa\, \Phi_{\beta,\tau,\gamma}(\kappa) \,|a_{\kappa}| < 1,\]
\noindent By using Remark \eqref{remark}, we give the estimate for the coefficients of an univalent function belong to $\mathcal{S}$ in $ \mathbb{U}$ also, by employ this estimate, we can get another  estimate for $ \ell_{1} $  in $ \mathcal{S}$ as follows,
 \begin{equation}\label{Eq 1}
\ell_{1} = \sum_{\kappa=2}^{\infty}\kappa \Phi_{\beta,\tau,\gamma} |a_{\kappa}| \leq \sum_{\kappa=2}^{\infty} \kappa^{2} \Phi_{\beta,\tau,\gamma}(\kappa) = \sum_{\kappa=2}^{\infty} \frac{(\kappa)^{2}}{\kappa!} (\Phi_{\beta,\tau,\gamma}(\kappa) \kappa!)=
 \sum_{\kappa=2}^{\infty} \frac{(\kappa)^{2}}{\kappa!} \ell(\kappa)< 1
 \end{equation}
 where
 \begin{equation}
 \ell(\kappa)= \frac{\Gamma(\frac{\beta}{\gamma+1}+1-\beta+\tau)}{ \Gamma(\frac{\beta}{\gamma+1}+1) }   \frac{\Gamma(\frac{\kappa+\beta-1}{\gamma+1}+1) ~(1)_{\kappa} } {\Gamma(\frac{\kappa+\beta-1}{\gamma+1}+1-\beta+\tau)}
 \end{equation}
 The series in \eqref{Eq 1} is transformed into a sum of two terms by employing the following relation:
   \begin{eqnarray}
 \frac{\kappa^{2}}{(1)_{\kappa}}= \frac{\kappa}{(1)_{\kappa-1}}= \frac{1}{(1)_{\kappa-1}} + \frac{1}{(1)_{\kappa-2}}
 \end{eqnarray}

 Depending on $ (1)_{\kappa}=\kappa! $ and $ (1)_{\kappa-1}=(\kappa-1)! $,  the estimate \eqref{Eq 1} becomes the next form:
 \begin{eqnarray*}
 \ell_{1} &\leq & \sum_{\kappa=2}^{\infty} \frac{\kappa^{2}}{(1)_{\kappa}} \ell(\kappa)
=\sum_{\kappa=2}^{\infty}\left(  \frac{1}{(1)_{\kappa-1}} + \frac{1}{(1)_{\kappa-2}}\right) \ell(\kappa) \\[0.6em]
&=& \sum_{\kappa=2}^{\infty}  \frac{\ell(\kappa)}{(1)_{\kappa-1}} + \frac{\ell(\kappa)}{(1)_{\kappa-2}}\\[0.6em]
&=& \sum_{\kappa=2}^{\infty} \frac{\Gamma(\frac{\beta}{\gamma+1}+1-\beta+\tau)}{ \Gamma(\frac{\beta}{\gamma+1}+1)}   \frac{\Gamma(\frac{\kappa+\beta-1}{\gamma+1}+1) } {\Gamma(\frac{\kappa+\beta-1}{\gamma+1}+1-\beta+\tau)} \frac{(1)_{\kappa}}{(1)_{\kappa-1}}
\\[0.6em]
&\qquad +& \sum_{\kappa=2}^{\infty} \frac{\Gamma(\frac{\beta}{\gamma+1}+1-\beta+\tau)}{ \Gamma(\frac{\beta}{\gamma+1}+1)}   \frac{\Gamma(\frac{\kappa+\beta-1}{\gamma+1}+1) } {\Gamma(\frac{\kappa+\beta-1}{\gamma+1}+1-\beta+\tau)} \frac{(1)_{\kappa}}{(1)_{\kappa-2}}
\\[0.6em]
&=& \frac{\Gamma(\frac{\beta}{\gamma+1}+1-\beta+\tau)}{ \Gamma(\frac{\beta}{\gamma+1}+1)} \bigg( \sum_{\kappa=1}^{\infty} \frac{\Gamma(\frac{\kappa+\beta}{\gamma+1}+1) } {\Gamma(\frac{\kappa+\beta}{\gamma+1}+1-\beta+\tau)} \frac{(1)_{\kappa+1}}{(1)_{\kappa}}\\[0.6em]
& \qquad +& \sum_{\kappa=0}^{\infty}
 \frac{\Gamma(\frac{\kappa+\beta+1}{\gamma+1}+1) } {\Gamma(\frac{\kappa+\beta+1}{\gamma+1}+1-\beta+\tau)} \frac{(1)_{\kappa+2}}{(1)_{\kappa}}\bigg)\\[0.6em]
  &= &\frac{\Gamma(\frac{\beta}{\gamma+1}+1-\beta+\tau)}{ \Gamma(\frac{\beta}{\gamma+1}+1)} \bigg( \sum_{\kappa=1}^{\infty} \frac{ \Gamma(\kappa+2) \Gamma(\frac{\kappa+\beta}{\gamma}+1) } {\Gamma(\frac{\kappa+\beta}{\gamma+1}+1-\beta+\tau)} \frac{1}{(1)_{\kappa}}\\[0.6em]
   &\qquad +& \sum_{\kappa=0}^{\infty}
 \frac{\Gamma(\kappa+3)\Gamma(\frac{\kappa+\beta+1}{\gamma+1}+1) } {\Gamma(\frac{\kappa+\beta+1}{\gamma+1}+1-\beta+\tau)} \frac{1}{(1)_{\kappa}}\bigg),
 \end{eqnarray*}
then by employing the Fox-Wright function given by \eqref{Fox}, we can transform the estimate $ \ell_{1}$ at $ z =1$ as follows
   \begin{eqnarray*}
={}_{2}\Psi_{1}\begin{scriptsize}
 \begin{bmatrix}
  (3,1), (1+\frac{\beta}{\gamma+1 }+\frac{1}{\gamma+1} ,\frac{1}{\gamma +1});   \\
\qquad  \qquad \qquad \qquad \qquad \qquad  1 \\
       (1-\beta +\tau+\frac{\beta}{\gamma +1} +       \frac{1}{\gamma+1},\frac{1}{\gamma +1});
           \end{bmatrix} \end{scriptsize}
  +
  {}_{2}\Psi_{1}  \begin{scriptsize}
    \begin{bmatrix}
  (2,1), (1+\frac{\beta}{\gamma+1 } ,\frac{1}{\gamma +1});            \\
 \qquad  \qquad \qquad \qquad \qquad 1\\
       (1-\beta +\upsilon+\frac{\beta}{\gamma +1} ,\frac{1}{\gamma +1});
           \end{bmatrix} \end{scriptsize}
           - \frac{ \Gamma(\frac{\beta}{\gamma+1}+1)} {\Gamma(\frac{\beta}{\gamma+1}+1-\beta+\tau)}    \\
 < \frac{ \Gamma(\frac{\beta}{\gamma+1}+1)}{\Gamma(\frac{\beta}{\gamma+1}+1-\beta+\tau)}
         \end{eqnarray*}
Hence,
\[ \Theta^{\beta,\tau,\gamma}: \mathcal{S} \rightarrow \mathcal{S}.\]
\end{proof}
 }\ \end{theorem}

 \begin{theorem}\rm{
\noindent Let the condition 1 as the Theorem \ref{Theo 1} is satisfied. If $ 0 \leq \beta-\tau<1$,
     \begin{eqnarray*}
  {}_{2}\Psi_{1} \begin{scriptsize}
  \begin{bmatrix}
  (2,1), (1+\frac{\beta}{\gamma+1 } ,\frac{1}{\gamma +1});            \\
  \qquad \qquad \qquad \qquad \qquad 1\\
       (1-\beta +\tau+\frac{\beta}{\gamma +1} ,\frac{1}{\gamma +1});           \\
\end{bmatrix}
  \end{scriptsize} <\,  2 \, \bigg( \frac{ \Gamma(\frac{\beta}{\gamma+1}+1)}{\Gamma(\frac{\beta}{\gamma+1}+1-\beta+\tau)}\bigg).
  \end{eqnarray*}
     then the operator  maps a convex function $ f(z)$ into a univalent function that is $ \Theta^{\beta,\tau,\gamma}: \mathcal{K}\rightarrow \mathcal{S}$.
      \begin{proof}
   Assume that $ f(z) \in \mathcal{K}$ and let
   \begin{eqnarray*}
 \Theta^{\beta,\tau,\gamma}f(z) =  z +
    \sum_{\kappa=2}^{\infty} w_{\kappa}
  z^{\kappa}
 \end{eqnarray*}
 be given by \eqref{extend}
 where
  \begin{eqnarray*}
  w_{\kappa}:=\Phi_{\beta,\tau,\gamma}(\kappa) a_{\kappa}
  \end{eqnarray*}
   and the function $\Phi_{\beta,\tau,\gamma}$ is given by \eqref{Psi}. To proof that the operator $\Theta^{\beta,\tau,\gamma}_{z}$ preserves the class $ \mathcal{S}$ we require the following sufficient condition
 \[\ell_{2};= \sum_{\kappa=2}^{\infty} \kappa |w_{\kappa}| = \sum_{\kappa=2}^{\infty} \kappa \Phi_{\beta,\tau,\gamma}(\kappa) |a_{\kappa}| < 1.\]
  We know That the coefficient  of a convex function belong to $\mathcal{S}$ is $ |a_{\kappa}|<1$. So  we can estimate $ \ell_{2} $  as follows,
 \begin{eqnarray}\label{Eq 2}
 \ell_{2} &=& \sum_{\kappa=2}^{\infty} \kappa \Phi_{\beta,\tau,\gamma} |a_{\kappa}| \leq \sum_{\kappa=2}^{\infty} \kappa^{2} \Phi_{\beta,\tau,\gamma}(\kappa)\\[0.6em]
 &= &\sum_{\kappa=2}^{\infty} \frac{(\kappa)^{2}}{\kappa!} (\Phi_{\beta,\tau,\gamma}(\kappa) \kappa!)=
 \sum_{\kappa=2}^{\infty} \frac{(\kappa)^{2}}{\kappa!} \ell(\kappa)< 1
  \end{eqnarray}
 where
\[
 \ell(\kappa)= \frac{\Gamma(\frac{\beta}{\gamma+1}+1-\beta+\tau)}{ \Gamma(\frac{\beta}{\gamma+1}+1) }   \frac{\Gamma(\frac{\kappa+\beta-1}{\gamma+1}+1) (1)_{\kappa} } {\Gamma(\frac{\kappa+\beta-1}{\gamma+1}+1-\beta+\tau)}
 \]
  and $ (a) _{\kappa} $ is Pochhammer symbol, with the following relation
\[
 \frac{\kappa}{(1)_{\kappa}}= \frac{1}{(1)_{\kappa-1}}.
 \]
  Since $ (1)_{\kappa}=\kappa ! $, then the estimate \eqref{Eq 2} has the next form
 \begin{eqnarray*}
 \ell_{2}
 &\leq & \sum_{\kappa=2}^{\infty} \frac{\kappa}{(1)_{\kappa}} \ell(\kappa)
= \sum_{\kappa=2}^{\infty}  \frac{1}{(1)_{\kappa-1}} \ell(\kappa)
 = \sum_{\kappa=2}^{\infty}  \frac{\ell(\kappa)}{(1)_{\kappa-1}}\\[0.6em]
&=& \sum_{\kappa=2}^{\infty} \frac{\Gamma(\frac{\beta}{\gamma+1}+1-\beta+\tau)}{\Gamma(\frac{\beta}{\beta+1}+1)}   \frac{\Gamma(\frac{\kappa+\beta-1}{\gamma+1}+1) } {\Gamma(\frac{\kappa+\beta-1}{\gamma+1}+1-\beta+\tau)} \frac{(1)_{\kappa}}{(1)_{\kappa-1}}
\\[0.6em]
&=& \frac{\Gamma(\frac{\beta}{\gamma+1}+1-\beta+\tau)}{ \Gamma(\frac{\beta}{\gamma+1}+1)}  \sum_{\kappa=1}^{\infty} \frac{\Gamma(\frac{\kappa+\beta}{\gamma+1}+1) } {\Gamma(\frac{\kappa+\beta}{\gamma+1}+1-\beta+\tau)} \frac{(1)_{\kappa+1}}{(1)_{\kappa}}
 \end{eqnarray*}
 which equivalents
 \begin{eqnarray*}
 = \frac{\Gamma(\frac{\beta}{\gamma+1}+1-\beta+\tau)}{ \Gamma(\frac{\beta}{\gamma+1}+1)} \sum_{\kappa=1}^{\infty} \frac{ \Gamma(\kappa+2) \Gamma(\frac{\kappa+\beta}{\gamma+1}+1) } {\Gamma(\frac{\kappa+\beta}{\gamma+1}+1-\beta+\tau)} \frac{1}{(1)_{\kappa}}.
 \end{eqnarray*}
Therefore, by utilizing the Fox-Wright function, we  transform the estimate $ \ell_{1}$ at $ z=1$,
\begin{eqnarray}
= \frac{\Gamma(\frac{\beta}{\gamma+1}+1-\beta+\tau)}{ \Gamma(\frac{\beta}{\gamma+1}+1)}  {}_{2}\Psi_{1} \begin{scriptsize}
\begin{bmatrix}
  (2,1), (1+\frac{\beta}{\gamma+1 } ,\frac{1}{\gamma +1});            \\
  \qquad \qquad \qquad \qquad \qquad 1\\
       (1-\beta+\tau+\frac{\beta}{\gamma +1} ,\frac{1}{\gamma +1});
           \end{bmatrix} \end{scriptsize} - 1  < 1.
   \end{eqnarray}
Hence
 \[\Theta^{\beta,\tau,\gamma}: \mathcal{K}\rightarrow \mathcal{S}.\]
 By this the proof is completed.
  \end{proof}}\end{theorem}

 \section{Boundedness and compactness in Bloch space}
\noindent In this section we characterize the boundedness and compactness of operator $\Theta^{\beta,\tau,\gamma}$  given by \eqref{normalization} on weighted $\mu$-Bloch space $ \mathbb{B}
^{\mu}_{w}$. 
First, let recall the well known Bloch space $ \mathbb{B}$ and weighted Bloch space $ \mathbb{B}_{w}$ \cite{z22,z23}
are defined  respectively,
\begin{definition}\label{Blo}\rm{ A holomorphic function $ f \in \mathcal{H}(\mathbb{U})$ is said to be in Bloch space $ \mathbb{B}$ whenever
\[
\| f \|_{\mathbb{B}}= \sup _{z \in   \mathbb{U}}\left( 1-\mid z \mid ^{2}\right)\mid f^{\prime}(z)\mid \, < \, \infty.
\]

and the little Bloch space $\mathbb{B}_{0}$ is given as follows
\[
 \lim_{|z|\rightarrow1-} \left( 1- |z|^{2}\right) |f^{\prime}(z)|=0.
\]
}\end{definition}
\begin{definition}\rm{
 Let $w: [0.1)\rightarrow [0,\infty)$ and $f$ be an analytic function on unit disk $ \mathbb{U}$ is said to be in  the weighted Bloch
 space $ \mathbb{B}_{w}$  if
\[
\left( 1- |z|\right)\mid f^{\prime}(z)\mid \, <  \hbar\, w(1-|z|), \quad z \in \mathbb{U}.
\]
for some $\hbar=\hbar_{f}\, >\, 0$. Not that, if $ w=1$ then $ \mathbb{B}_{w}$ reduces to the classical Bloch space $\mathbb{B}$.
 Further, the weighted $ \mu$- Bloch space $ \mathbb{B}_{w}^{\mu}$, covering of all $ f\in \mathbb{B}_{w}^{\mu}$ defined by
\[\lim_{|z|\rightarrow -1}\frac{(1-|z|)^{\mu} |f^{\prime}(z)|}{w(1-|z|)}=0.  \]
and
\begin{eqnarray}\label{wBloch}
\| f \|_{\mathbb{B}_{w}^{\mu}}= \sup _{z \in   \mathbb{U}} |f^{\prime}(z)|  \frac{\left( 1- | z |\right)^{\mu}}{w(1-|z|)}\, < \, \infty.
\end{eqnarray}
It is easy to note that if an analytic function $ g(z)\in \mathbb{B}_{w}^{\mu}$, then
\begin{eqnarray}\label{inequality}
\sup_{z \in \mathbb{U}}|\Bbbk  g(z)| \frac{(1-|z|)^{\mu}}{w(1-|z|)} \leq c\, < \infty.
\end{eqnarray}
where $\Bbbk$ is a \textit{positive} number.
}\end{definition}

\begin{theorem}\label{Th bounded}\rm{
Let $f$ be an analytic function on open unit disk $ \mathbb{U}$, and $ \mathbb{B}_{w}^{\mu}; \,w : [0,1) \rightarrow [0,\infty)$. Then
\[ f \in \mathbb{B}_{w}^{\mu} \Longleftrightarrow \Theta^{\beta,\tau,\gamma} f \in \mathbb{B}_{w}^{\mu}.\]
\begin{proof}

Let suppose $ f \in \mathbb{B}_{w}^{\mu}$, then by using equalities \eqref{normalization} and \eqref{wBloch}, we obtain
\begin{eqnarray*}
|| \Theta^{\beta,\tau,\gamma}f||_{\mathbb{B}_{w}^{\mu}}
  &=&\sup_{z \in \mathbb{U}}\bigg|\left(z\left[\frac{ \Gamma(\frac{1+\beta}{\gamma+1}+1)}{ \Gamma(\frac{1+\beta}{\gamma+1}+1-\beta+\tau)}{}_{2}\Psi_{1}[z] * f(z)\right]\right)^{\prime} \bigg| \frac{(1-|z|)^{\mu}}{w(1-|z|)}
  \\[0.6em]
 &\leq & \sup_{z \in \mathbb{U}}\bigg|\tfrac{ \Gamma(\frac{1+\beta}{\gamma+1}+1)}{ \Gamma(\frac{1+\beta}{\gamma+1}+1-\beta+\tau)}{}_{2}\Psi_{1}[z] * f(z) \bigg| \frac{(1-|z|)^{\mu}}{w(1-|z|)}
  \\[0.6em]
 &\qquad +& \sup_{z \in \mathbb{U}}\bigg| z \bigg(\tfrac{ \Gamma(\frac{1+\beta}{\gamma+1}+1)}{ \Gamma(\frac{1+\beta}{\gamma+1}
+ 1-\beta+\tau)}{}_{2}\Psi_{1}[z] * f(z)\bigg)^{\prime} \bigg| \frac{(1-|z|)^{\mu}}{w(1-|z|)},
  \end{eqnarray*}
  by following the derivative of convolution due to Ruscheweyh (see \cite{z44}, pp. 39), we have
 \begin{eqnarray*}
|| \Theta^{\beta,\tau,\gamma}f||_{\mathbb{B}_{w}^{\mu}}
&\leq &
 \sup_{z \in \mathbb{U}}
 \bigg|\frac{ \Theta^{\beta,\tau,\gamma}f}{z}\bigg| \frac{(1-|z|)^{\mu}}{w(1-|z|)}
 +  \sup_{z \in \mathbb{U}}
 \bigg| \tfrac{ \Gamma(\frac{1+\beta}{\gamma+1}+1)}{ \Gamma(\frac{1+\beta}{\gamma+1}
+ 1-\beta+\tau)}{}_{2}\Psi_{1}[z] *z f^{\prime}(z) \bigg| \frac{(1-|z|)^{\mu}}{w(1-|z|)},
\\[0.6em]
&\leq &
 c
 +  c \,\sup_{z \in \mathbb{U}}
 \big| f^{\prime}(z) \big| \frac{(1-|z|)^{\mu}}{w(1-|z|)},\quad |z|<1
\end{eqnarray*}
where $ {}_{2}\Psi_{1}$ is given by \eqref{q: G} and $ |z| <1$. Hence $ \Theta^{\beta,\tau,\gamma}f \in \mathbb{B}_{w}^{\mu} $. On the other hand,
if $ \Theta^{\beta,\tau,\gamma}f \in \mathbb{B}_{w}^{\mu} $, then
\[ || \Theta^{\beta,\tau,\gamma}f||_{  \mathbb{B}_{w}^{\mu}}= \sup_{z \in \mathbb{U}}|(\Theta^{\beta,\tau,\gamma}f)^{\prime}(z)| \frac{(1-|z|)^{\mu}}{w(1-|z|)} < \infty.
\]
Since \[ \Theta^{\beta,\tau,\gamma}f(z) =  \frac{\Gamma(\frac{1+\beta}{\gamma+1}+1)}{\Gamma(\frac{1+\beta}{\gamma+1}+1-\beta+\tau)}z \, {}_{2}\Psi_{1}[z] * f(z),
\]
then,
\[ f(z)= \frac{\Gamma(\frac{1+\beta}{\gamma+1}+1-\beta+\tau)}{ \Gamma(\frac{1+\beta}{\gamma+1}+1)z \,{}_{2}\Psi_{1}[z]} \Theta^{\beta,\tau,\gamma}f(z)
\]

is an analytic function on open unit disk $ \mathbb{U}$,
 we have the following
\[ \sup_{z \in \mathbb{U}_{\frac{1}{2}}} |f^{\prime}(z)| \frac{ (1-|z|)^{\mu}}{w(1-|z|)} \leq c_{1}< \infty.
\]
where $ \mathbb{U}_{\frac{1}{2}} = \lbrace z: |z|<  \frac{1}{2} \rbrace$ and $ c_{1}$ is a \textit{positive} constant. Therefore by $ \Theta^{\beta,\tau,\gamma}f(z) \in \mathbb{B}_{w}^{\mu}$, we obtain
\begin{small}
\begin{eqnarray*}
|| f||_{\mathbb{B}_{w}^{\mu}} &=& \sup_{z \in \mathbb{U}} |f^{\prime}(z)| \frac{ (1-|z|)^{\mu}}{w(1-|z|)}
 \\[0.6em]
 &\leq &  c_{1} +  \sup_{z \in \mathbb{U}\setminus \mathbb{U}_{\frac{1}{2}}} \bigg|\left(\tfrac{\Gamma(\frac{1+\beta}{\gamma+1}+1-\beta+\tau) } { \Gamma(\frac{1+\beta}{\gamma+1}+1)} \frac{\Theta^{\beta,\tau,\gamma}
 f}{z\, {}_{2}\Psi_{1}[z]} \right)^{\prime}\bigg| \frac{ (1-|z|)^{\mu}}{w(1-|z|)},
 \\[0.6em]
 &=&c_{1}+ \sup_{z \in \mathbb{U}\setminus \mathbb{U}_{\frac{1}{2}}} \bigg|\tfrac{\Gamma(\frac{\beta}{\gamma+1}+1-\beta+\tau)}{\Gamma(\frac{\beta}{\gamma+1}+1)} \left\lbrace \frac{z\, {}_{2}\Psi_{1}[z](\Theta^{\beta,\tau,\gamma}f)^{\prime}-\Theta^{\beta,\tau,\gamma}f(z\, {}_{2}\Psi_{1}[z])^{\prime}}{ z^{2}\,({}_{2}\Psi_{1}[z])^{2}} \right\rbrace
 \bigg| \frac{ (1-|z|)^{\mu}}{w(1-|z|)},
 \\[0.6em]
   & \leq &  c_1 +\frac{ 2}{c} \sup_{z \in \mathbb{U}\setminus \mathbb{U}_{\frac{1}{2}}} | \Theta^{\beta,\tau,\gamma}
 f^{\prime}| \frac{ (1-|z|)^{\mu}}{w(1-|z|)} + \frac{2}{ c}\sup_{z \in \mathbb{U}\setminus \mathbb{U}_{\frac{1}{2}}} |  f * {}_{2}\Psi_{1}[z]| \frac{ (1-|z|)^{\mu}}{w(1-|z|)}
  \\[0.6em]
 & \leq &  c_1 + 2|| \Theta^{\beta,\tau,\gamma}
 f||_{B_w}^{\mu}+  2 c_2   < \infty.
\end{eqnarray*}
\end{small}
where $ c_{1}$ and $ \sup_{z \in \mathbb{U}\setminus \mathbb{U}_{\frac{1}{2}}} \frac{ |  f * {}_{2}\Psi_{1}[z]|}{c} \frac{ (1-|z|)^{\mu}}{w(1-|z|)}  \leq c_2 $ are \textit{positive} constant, so $ f \in \mathbb{B}_{w}^{\mu}$.
This completes the proof.
\end{proof}
}\end{theorem}
 \begin{theorem}\rm{
 let $f$ be an analytic function on open unit disk $ \mathbb{U}$, and $ \mathbb{B}_{w}^{\mu}; \,w : [0,1) \rightarrow [0,\infty)$. Then \[ \Theta^{\beta,\tau,\gamma} f:\mathbb{B}_{w}\rightarrow\mathbb{B}_{w}^{\mu}
 \]
 is compact.
 \begin{proof}
If $ \Theta^{\beta,\tau,\gamma}f $ is compact, then it is bounded and by Theorem \ref{Th bounded} it satisfies that
 $ f \in  \mathbb{B}_{w}$ because $ \mathbb{B}_{w} \subset \mathbb{B}_{w}^{\mu}$.
  Let assume that $ f \in \mathbb{B}_{w}$, that $( f_n)_{n\in \mathbb{N}} \subset \mathbb{B}_{w}^{\mu}$ be such that $ f_{n}\rightarrow 0$ converges uniformly on  $\overline{ \mathbb{U}}$ as $ n \rightarrow   \infty$. Since $ (f)_{n \in \mathbb{N}}$ convergence uniformly on each compact  $ \mathbb{U}$, we have that there in $ \mathbb{N}>0$ such that for every $ n > \mathbb{N}$ and every $ z \in \mathbb{U}$, there is an $ 0<\delta <1$, such that for every $ n \geq 1$,
 \[\bigg|\tfrac{\Gamma(\frac{\beta}{\gamma+1}+1-\beta+\tau)}{\Gamma(\frac{\beta}{\gamma+1}+1)}
   {}_{2}\Psi_{1}[z]*f_{n}(z)
  \bigg| < \varepsilon
 \]
 where $ \delta < |z|< 1$, Since $ \delta$ is arbitrary, then we can choose
  \[ \big|\frac{(1-|z|)}{(1-|z|)^{\mu}}\big|<1
  \]
 for all $ \delta < |z|<1$ and
 \begin{eqnarray}
 ||\Theta^{\beta,\tau,\gamma} f_{n}||_{\mathbb{B}_{w}}&=& \sup_{z \in \mathbb{U}\setminus \mathbb{U}_{\delta}} \big|\tfrac{\Gamma(\frac{\beta}{\gamma+1}+1-\beta+\tau)}{\Gamma(\frac{\beta}{\gamma+1}+1)}
 \lbrace
 \big(z[ {}_{2}\Psi_{1}[z]*f_{n}(z)]\big)^{\prime}
 \rbrace
 \big| \frac{ (1-|z|)}{w(1-|z|)},\nonumber \\[0.6em]
 &= &
 \sup_{z \in \mathbb{U}\setminus \mathbb{U}_{\delta}} \big|\tfrac{\Gamma(\frac{\beta}{\gamma+1}+1-\beta+\tau)}{\Gamma(\frac{\beta}{\gamma+1}+1)}
   {}_{2}\Psi_{1}[z]*f_{n}(z)
  \big| \frac{ (1-|z|)}{w(1-|z|)} \nonumber
 \\[0.6em]
 &\qquad  + &  \sup_{z \in \mathbb{U}\setminus \mathbb{U}_{\delta}} \big|\tfrac{\Gamma(\frac{\beta}{\gamma+1}+1-\beta+\tau)}{\Gamma(\frac{\beta}{\gamma+1}+1)}
   {}_{2}\Psi_{1}[z]* z f^{\prime}_{n}(z)
  \big| \frac{ (1-|z|)}{w(1-|z|)} \nonumber
 \\[0.6em]
 & \leq & \varepsilon +
   \sup_{z \in \mathbb{U}\setminus \mathbb{U}_{\delta}} \big|\tfrac{\Gamma(\frac{\beta}{\gamma+1}+1-\beta+\tau)}{\Gamma(\frac{\beta}{\gamma+1}+1)}
   {}_{2}\Psi_{1}[z]* z f^{\prime}_{n}(z)
  \big| \frac{ (1-|z|)^{\mu}}{w(1-|z|)},\quad |z|<1 \nonumber
 \\[0.6em]
 &  \leq & \varepsilon +
   c ||f_{n}||_{\mathbb{B}_{w}^{\mu} }\label{last}
  \end{eqnarray}
Since for $ f_{n}\rightarrow 0$ on $\overline{ \mathbb{U} }$ we get $  ||f_{n}||_{\mathbb{B}_{w}^{\mu} }\rightarrow 0$, and that $ \varepsilon$ is an arbitrary positive number, by letting $ n\rightarrow \infty$ in
 \eqref{last}, we have that $ \lim_{n\rightarrow\infty}||\Theta^{\beta,\tau,\gamma}f_{n}||_{\mathbb{B}_{w}}=0$. Thus $ \Theta^{\beta,\tau,\gamma}$ is compact.
 \end{proof}
 }\end{theorem}

\end{document}